\documentclass[a4paper,11pt]{article}

\usepackage[english]{babel}
\usepackage[english]{translator}
\usepackage[T1]{fontenc}
\usepackage[utf8]{inputenc}

\usepackage{amsmath}
\usepackage{amssymb}
\usepackage{amsfonts}
\usepackage{amstext}
\usepackage{amsthm}
\usepackage{bbm}

\usepackage{graphicx}
\usepackage{color}
\usepackage{psfrag}
\usepackage{subfigure}
\usepackage{float}

\usepackage{marvosym}
\usepackage{enumerate}
\usepackage{enumitem}
\usepackage{todonotes}
\definecolor{myForestGreen}{RGB}{34,139,34}
\newcommand{\corr}[1]{{\color{red} #1}} 

\usepackage{hyperref}

\usepackage[
  hmarginratio={1:1},     
  vmarginratio={1:1},     
  textwidth=460pt,        
  heightrounded,          
  bottom=3.65cm,
  top=3.65cm,
]{geometry}

\hypersetup{
  pdffitwindow=false,
  pdfhighlight=/O,
  pdfnewwindow,
  colorlinks=true,
  citecolor=red,             
  linkcolor=blue,            
  menucolor=blue,            
  urlcolor=blue,             
  pdfpagemode=UseOutlines,
  bookmarksnumbered=true,
  linktocpage,
  pdfkeywords={},
  pdfcreator={pdflatex},
  pdfproducer={LaTeX with hyperref}
}

\newcommand{\D}{\mathcal{D}}
\newcommand{\C}{\mathbb{C}}
\newcommand{\R}{\mathbb{R}}
\newcommand{\N}{\mathbb{N}}
\renewcommand{\S}{\mathcal{S}}

\newcommand{\T}{\mathcal{T}}

\renewcommand{\ker}{\mathrm{ker}}

\renewcommand{\i}{\mathrm{i}}

\renewcommand{\H}{\mathcal{H}}

\newcommand{\ls}{\lesssim}
\newcommand{\LOD}{{\text{\tiny LOD}}}

\newcommand{\Rlod}{R_{\LOD}}

\renewcommand{\P}{\mathbb{P}}
\renewcommand{\corr}{\mathcal{C}}
\newcommand{\gs}{\mathrm{gs}}

\newcommand{\sect}
{
  \setcounter{equation}{0}
  \setcounter{figure}{0}
  \section
}

\theoremstyle{definition}
\newtheorem{definition}{Definition}[section]

\newtheorem{remark}[definition]{Remark}

\theoremstyle{plain}

\allowdisplaybreaks

\begin{document}
\begin{center}
{\LARGE 
Localized Orthogonal Decomposition Methods vs. Classical FEM for the Gross-Pitaevskii Equation
}
\end{center}

\begin{center}
{\large Christian D\"oding\footnote[1]{Institute for Numerical Simulation, University of Bonn, D-53115 Bonn, Germany, \\ e-mail: \textcolor{blue}{doeding@ins.uni-bonn.de}.}}\\[2em]
\end{center}

\noindent
\begin{center}
\begin{minipage}{0.8\textwidth}
  {\small
    \textbf{Abstract.} The time-dependent Gross-Pitaevksii equation (GPE) is a nonlinear Schr\"odinger equation which is used in quantum physics to model the dynamics of Bose-Einstein condensates. In this work we consider numerical approximations of the GPE based on a multiscale approach known as the localized orthogonal decomposition. Combined with an energy preserving time integrator one derives a method which is  of high order in space and time under mild regularity assumptions. In previous work, the method has been shown to be numerically very efficient compared to first order Lagrange FEM. In this paper, we further investigate the performance of the method and compare it with higher order Lagrange FEM. For rough problems we observe that the novel method performs very efficient and retains its high order, while the classical methods can only compete well for smooth problems.
}
\end{minipage}
\end{center}



 
\sect{Introduction}
\label{sec:intro}
This work is devoted to numerical approximations of the time-dependent Gross-Pitaevskii equation (GPE) on a bounded convex polyhedral domain $\D \subset \R^d$, $d = 1,2,3$ that seeks for a complex-valued function $u = u(x,t) \in \C$ such that
\begin{align*}
	\i \partial_t u = - \Delta u + Vu + \beta |u|^2 u
\end{align*}
holds in $\D \times [0,T]$ for some $T > 0$, $\beta \in \R$, and a real-valued positive function $V = V(x) \ge 0$. For the well-posedness of the problem we assume homogeneous Dirichlet boundary conditions, i.e., $u(\cdot,t) = 0$ on $\partial \D$, and a sufficiently smooth initial value $u(\cdot,0) = u_0$. The GPE is used in applications to model the dynamics of so-called Bose-Einstein condensates (BECs) and plays an important role in modern physics, see e.g. \cite{BaoCai} and the references therein. \\
In this paper we consider numerical approximations of the GPE based on a spatial discretization by localized orthogonal decomposition (LOD, cf. \cite{MalqvistPeterseim14}) combined with an energy preserving high order time integrator from \cite{Makridakis}. We compare the considered method, first introduced in \cite{DoedingHenningWarnegard}, with classical spatial discretizations based on standard higher order Lagrange FEM (cf. \cite{Makridakis,DoedingHenning}). In preliminary work (cf. \cite{DoedingHenningWarnegard,HenningWarnegard}) a spatial discratization by LOD turned out to be very efficient in solving the time-dependent GPE. This is due to the fact that the spatial discretization by LOD gives a third order method in space w.r.t. the energy norm under mild regularity assumptions. Therefore, a combination with an energy preserving high order time integrator is likely to give highly accurate approximations. Nevertheless, the method has usually been compared to first order Lagrange FEM (see \cite{DoedingHenningWarnegard,HenningWarnegard}) and one might think that higher order Lagrange FEM would perform similarly as the considered LOD method for smooth problems. \\
The contribution of this work is a comparison of the aforementioned methods in numerical experiments. We demonstrate that the classical methods can compete well with the proposed LOD method when solving smooth problems. However, for rough problems with low regularity, we show that the third order convergence of the LOD method is preserved whereas the classical methods obviously fail due to the lack of regularity. In this case, the desired LOD method can play to its strengths and perform much better than classical methods, making it an efficient choice as a solver for the GPE and the simulation of BECs.  \\ 
In the following sections, we first outline the analytical framework in which we are working and recall the approach of the LOD, before deriving the fully-discrete system through a combination with a continuous Galerkin time integrator. We conclude the preliminaries by giving the corresponding fully-discrete methods based on higher order Lagrange FEM, which we compare in the numerical experiments. In Section \ref{sec:comp}, we investigate the performance of the presented methods for different problems and discuss the numerical observations.

\sect{Problem Setting}
\label{sec:setting}

We denote by $L^p(\D) = L^p(\D,\C)$, $1 \le p \le \infty$ the standard Lebesgue spaces of complex-valued $p$-integrable functions and equip the complex Hilbert space $L^2(\D)$ with the inner product $(u,v) = \int_\D \overline{u} \, v \, dx$ which is conjugate-linear in its first argument. Further, let $H^1_0(\D)$ denote the usual Sobolev space of weakly differentiable functions in $L^2(\D)$ with zero trace and let $H^{-1}(\D)$ be its dual space. \\
The variational form of the GPE seeks a function $u \in L^\infty([0,T],H^1_0(\D))$ with $\partial_t u \in L^\infty([0,T], H^{-1}(\D))$ such that $u(\cdot,0) = u_0$ and for almost every $t \in (0,T]$
\begin{align} \label{GPE}
	( \i \partial_t u, v) = (\nabla u, \nabla v) + (V u , v) + \beta (|u|^2 u, v) \quad \forall v \in H^1_0(\D).
\end{align}
If we assume from now on that $V \in L^\infty(\D,[0,\infty))$, $\beta \ge 0$, and $u_0 \in H^1_0(\D)$ then \eqref{GPE} admits a unique solution for some $T > 0$ (depending on $u_0$) and we refer to \cite{Cazenave} for further details about the well-posedness of \eqref{GPE}. In particular, we note that such a solution is continuous in time w.r.t. $L^2(\D)$ so that the initial condition makes sense. The energy of the system \eqref{GPE} is given by the so-called Gross-Pitaevskii energy functional
\begin{align*}
	E(v) = \frac{1}{2} \int_\D |\nabla v|^2 + V |v|^2 + \frac{\beta}{2} |v|^4 dx, \quad v \in H^1_0(\D)
\end{align*}
which is known to be conserved by the time evolution of \eqref{GPE}, i.e., $E(u(t)) = E(u_0)$ for all $t \in [0,T]$. In fact \eqref{GPE} can be formally derived as a Hamiltonian system from the energy $E$, see \cite{Bao,BaoCai}, so that the energy conservation follows immediately from the derivation of the system. 
From the numerical point of view it is of utmost importance to ensure the conservation of energy by numerical approximations as accurate as possible (but staying economically). \\
We introduce the conjugate-linear operator
\begin{align*}
	\H : H^1_0(\D) \rightarrow H^{-1}(\D) \quad \langle \H v, w \rangle := (\nabla v, \nabla w) + (Vv, w)
\end{align*}
where $\langle \cdot, \cdot \rangle$ denotes the dual pairing between $H^{-1}(\D)$ and $H^1_0(\D)$. In particular, $\H$ induces the sesquilinear form
\begin{align*}
	(v,w)_{\H} := (\nabla v, \nabla w) + (Vv, w)
\end{align*}
given by the linear part of the GPE which clearly is bounded and elliptic on $H^1_0(\D)$ since $V$ is assumed to be bounded and positive. Hence, by Lax-Milgram lemma the inverse operator $\H^{-1}: H^{-1}(\D) \rightarrow H^1_0(\D)$ exists and is bounded so that $v = \H^{-1} f$ is the solution of the linear elliptic equation $(v,w)_\H = \langle f,w \rangle$, $w \in H^1_0(\D)$ for any $f \in H^{-1}(\D)$.

\sect{Spatial Discretization by Localized Orthogonal Decomposition}

In this section we briefly recall the basics of the spatial discretization by localized orthogonal decomposition in the context of the GPE. Further details can be found in \cite{AltmannHenningPeterseim,EngwerHenningMalqvistPeterseim,HenningMalqvist,MalqvistPeterseim14,MalqvistPeterseim21} for the LOD in general settings and in \cite{DoedingHenningWarnegard,HenningWarnegard} for its particular application to the GPE. \\
We start by choosing a family $\T_H$ of quasi-uniform and shape-regular triangulations of $\D$ with mesh size parameter $0 < H < 1$ and let
\begin{align*}
	\P^1(\T_H) := \{ v \in C(\overline{D}) \cap H^1_0(\D): v_{|K} \in \mathbb{P}^1(K), \, K \in \T_H \},
\end{align*}
be the family of $H^1$-conforming $\P^1$-Lagrange finite element spaces associated with $\T_H$. The ideal approximation space for the (localized) orthogonal decomposition is defined by
\begin{align} \label{LODspace}
	V_{\LOD} := \H^{-1} \P^1(\T_H).
\end{align}
By construction we have $\dim V_{\LOD} = \dim \P^1(\T_H) =: N_H$ and with the $L^2$-projection $P_H: L^2(\D) \rightarrow \P^1(\T_H)$, i.e., $(v - P_H v,w) = 0$ for all $w \in \P^1(\T_H)$, it is easy to show (cf. \cite{AltmannHenningPeterseim}) that
\begin{align*}
	V_{\LOD} \,\, \bot_{(\cdot,\cdot)_{\H}} \,\, W := \ker (P_H) \cap H^1_0(\D).
\end{align*}
Indeed, this new space comes with very nice approximation properties: Let $R_\LOD: H^1_0(\D) \rightarrow V_\LOD$ be the Ritz projection onto $V_\LOD$, i.e., $( v - R_\LOD v, w)_\H = 0$ for all $w \in V_\LOD$. Then it was shown, e.g. in \cite{HenningMalqvist,MalqvistPeterseim14}, that if $v \in H^1_0(\D)$ is such that $\H v \in H^1_0(\D) \cap H^2(\D)$ it satisfies the estimate
\begin{align} \label{LODest}
	\| v - \Rlod v \|_{H^1} \ls H^3.
\end{align} 
We see immediately, that if $v \in H^1_0(\D) \cap H^4(\D)$ (which is in general stronger than $\H v \in H^1_0(\D) \cap H^2(\D)$) the third order approximation in \eqref{LODest} is also obtained by a standard $\P^3$-Lagrange FEM. Therefore, one might expect that a $\P^3$-Lagrange FEM could compete well with the LOD in solving \eqref{GPE} if the problem is sufficiently smooth. However, on a fixed mesh, i.e., $H > 0$ fixed, the dimension of $V_\LOD$ is equal to the dimension of the $\P^1$-Lagrange FE space which is smaller than the dimension of the $\P^3$-Lagrange FE space on the same mesh. Hence, the time stepping using a LOD spatial discretization is expected to be faster due to the lower dimensionality of the system matrices. Nevertheless, as we will see later in Section \ref{sec:comp}, the assembly of the nonlinear terms in the LOD space may require additional computational costs in each time step so that the lower dimensionality of $V_\LOD$ may not pay off in the overall computation. \\
In contrast, for rough potentials $V \in L^\infty$ without any further regularity we cannot expect $H^4$-regularity of the solution $u$, but $H^2$-regularity of $\H u$. Therefore, a method based on LOD is expected to preserve the third order convergence in space due to \eqref{LODest}, whereas a $\P^3$-Lagrange FEM is unsuitable due to the lack of regularity. \\
A natural choice for a basis of $V_\LOD$ is $\H^{-1} \varphi_i^H$, $i = 1,\dots,N_H$ where $\varphi^H_i \in \P^1(\T_H)$ denotes the standard hat function associated with the $i$th node of $\T_H$. But this basis (and in general any other basis) is no longer locally supported, which, in contrast to its good approximation properties, gives the impression that $V_\LOD$ is not a suitable approximation space as it implies densely populated system matrices. However, it is known that $V_\LOD$ is spanned by a basis whose functions decay exponentially (cf. \cite{MalqvistPeterseim14,MalqvistPeterseim21}) and can therefore be restricted to local patches by making small localization errors. We make this more precise in the following: \\
We start by considering $V_\LOD$ from \eqref{LODspace} as a correction of $\P^1(\T_H)$. Let $\corr: H^1_0(\D) \rightarrow W$ be the corrector such that
\begin{align*}
	(\corr v, w)_\H = (v, w)_\H \quad \forall w \in W.
\end{align*}
Then one shows that $V_\LOD$ is equivalently characterized via
\begin{align*}
	V_\LOD = (I - \corr) \P^1(\T_H)
\end{align*}
and is spanned by the (exponentially decaying) basis $\varphi_i^\LOD := (I - \corr) \varphi^H_i$, $i = 1,\dots N_H$. For $\ell \in \N$ we introduce the $\ell$-patch of an element $K \in \T_H$ recursively by
\begin{align*}
	N^\ell(K) := \bigcup \{ \hat K \in \T_H: \hat K \cap N^{\ell -1}(K) \neq \emptyset \}, \quad N^0(K) := K
\end{align*}
and the subspaces
\begin{align*}
	W^\ell_K := \{ v \in W: \mathrm{supp}(v) \subset N^\ell(K) \} \subset W.
\end{align*}
Further, we define the local element corrector $\corr^\ell_K:H^1_0(\D) \rightarrow W^\ell_K$ of $K$ via
\begin{align} \label{localproblem}
	(\corr^\ell_K v, w)_{\H,N^\ell(K)} = (v,w)_{\H,K} \quad \forall w \in W^\ell_K
\end{align}
where $(v, w)_{\H,S} := \int_{S} \overline{\nabla v} \cdot \nabla w + V \overline{v} w \, dx$ for $\S \subset \D$. This leads us to the localized corrector
\begin{align*}
	\corr^\ell: H^1_0(\D) \rightarrow W, \quad \corr^\ell v = \sum_{K \in \T_H} \corr^\ell_K v
\end{align*}
and we define the approximation space of the localized orthogonal decomposition as $V^\ell_\LOD := (I - \corr^\ell) \P^1(\T_H)$ with basis functions $(I - \corr^\ell) \varphi^H_i$, $i = 1,\dots N_H$ which have local support with a diameter of the order $\mathcal{O}(\ell H)$. The localization error is moderate since it follows, e.g. from \cite{HenningMalqvist,HenningWarnegard}, that $\| \corr^\ell v - \corr v \|_{H^1} \ls \exp(-\rho \ell) \| v \|_{H^1}$ for some $\rho > 0$. Hence, if we choose $\ell > 0$ sufficiently large in our experiments (depending on $\tau$ and $H$) the error caused by the localization is negligible compared with the dominating error from \eqref{LODest}.
\begin{remark}
We emphasize that there are other localization strategies, such as the novel super-localization strategy developed in \cite{HauckPeterseim}. Recent experiments, see \cite{PeterseimWarnegardZimmer}, have shown that this approach also improves the classical LOD approach in the context of the GPE. This suggests significant margin for improvements in the results presented here. However, we do not focus on the effect of localization strategies in this paper and therefore decided to use the more established classical localization strategy for the sake of brevity.  
\end{remark}

For the implementation of the LOD one has to solve the local problems \eqref{localproblem} for $v = \varphi^H_i$ which is done in practice by a fine scale discretization of mesh size $0 < h < H$, see e.g. \cite{EngwerHenningMalqvistPeterseim}. We choose a standard $\P^1$-Lagrange FEM on a fine triangulation $\T_h$ which we assume for simplicity to be a refinement of the coarse mesh $\T_H$. This leads us to the fully-discrete spatial approximation space
\begin{align*}
	V^{\ell,h}_{\LOD} := (I - \corr^{\ell,h}) \P^1(\T_H)
\end{align*}
where the corrector $\corr^{\ell,h}: H^1_0(\D) \rightarrow W_h := W \cap \P^1(\T_h)$ is given by
\begin{align*}
	\corr^{\ell,h} v = \sum_{K \in \T_H} \corr^{\ell,h}_K v
\end{align*}
with
\begin{align*}
	(\corr^{\ell,h}_K v, w)_{\H,N^\ell(K)} = (v,w)_{\H,K} \quad \forall w \in W^\ell_K \cap \P^1(\T_h). 
\end{align*}

\sect{Time Integration and Fully-Discrete Schemes}
\label{subsec:time}

As a next step, we combine the spatial discretization by LOD with a suitable energy conserving time integrator to obtain a fully-discrete method for solving \eqref{GPE}. The time integrator is designed as a continuous Galerkin method in time which was first introduced in \cite{Makridakis} for the NLS (i.e. $V = 0$) in two space dimensions and combined with standard Lagrange FE. The method was further analyzed for the GPE in up to three space dimensions in \cite{DoedingHenning} and a combination with a LOD space discretization was first proposed in \cite{DoedingHenningWarnegard}. \\
For the time integrator we define time instances $t_n = n \tau$ for $n= 0,\dots,N$ with step size $\tau = T/N > 0$ and subintervals $I_n = (t_n,t_{n+1}]$. For a polynomial degree $q \in \N$ we define the spaces
\begin{align*}
	& \tilde{W}^{\ell,h}_{q,\LOD} := \left\{ v : [0,T] \rightarrow V_\LOD^{\ell,h} : v_{| I_n} \in \P^q(I_n), \, n = 0,\dots,N-1 \right\}, \\
	& W^{\ell,h}_{q,\LOD} := \tilde{W}^{\ell,h}_{q,\LOD} \cap C([0,T],V^{\ell,h}_{\LOD}),
\end{align*}
where $\P^q(I_n)$ denotes the space of formal polynomials of degree $q \in \N$ on $I_n$ with values in $V^{\ell,h}_\LOD$. The fully-discrete method, which we refer to as the cG-LOD-method, now seeks $u_\LOD \in W^{\ell,h}_{q,\LOD}$ with $u_\LOD(\cdot,0) = P_\LOD(u_0)$ and such that
\begin{align} \label{cG-LOD-scheme}
	\int_{0}^T (\i \partial_t u_{\LOD}, v) - (u_\LOD,v)_\H - \beta \big( P_{\LOD}( |u_\LOD |^2) u_\LOD, v \big) \, dt = 0 \quad \forall v \in \tilde{W}^{\ell,h}_{q-1,\LOD}.
\end{align}
Using the continuity in time this solution can be computed recursively on each $I_n$. We refer to \cite{DoedingHenningWarnegard} for a detailed presentation on how the time stepping scheme \eqref{cG-LOD-scheme} is realized in an implementation. 

\begin{remark}
We emphasize the projection $P_\LOD (|u_\LOD|^2)$ onto $V_\LOD$ of the density of the solution in the nonlinear term. As elaborated in \cite{DoedingHenningWarnegard,HenningWarnegard} this projection allows for a speed up in the performance of the scheme in contrast to the naive method where the nonlinear term in \eqref{cG-LOD-scheme} is replaced by $(|u_\LOD|^2u_\LOD, v)$. Projecting the density onto the LOD space allows for a more efficient assembly of the nonlinear terms during the time stepping via a pre-computed tensor of the form $\omega_{ijk} = (\varphi_i^\LOD \varphi_j^\LOD, \varphi_k^\LOD)$ where $\varphi_i^\LOD$ denotes the basis functions of $V_\LOD$, cf. \cite{DoedingHenningWarnegard}.
\end{remark}

By testing with $\partial_t u_\LOD \in \tilde{W}^{\ell,h}_{q-1,\LOD}$ one infers that the cG-LOD-method \eqref{cG-LOD-scheme} preserves the modified energy
\begin{align*}
	E_\LOD(v) = \frac{1}{2} \int_\D |\nabla v|^2 + V |v|^2 + \frac{\beta}{2} P_\LOD(|v|^2)|v|^2 dx, \quad v \in H^1_0(\D)
\end{align*}
at the discrete time instances $t_n$, i.e., $E_\LOD(u_\LOD(t_n)) = E_\LOD(P_\LOD(u_0))$ for all $n=0,\dots,N$. Indeed, preserving only the modified energy $E_\LOD$ instead of the original energy $E$ does not cause relevant errors since it was shown in \cite{HenningWarnegard} that $|E_{\LOD}(v) - E(v)| \ls H^6$ if $v$ and $V$ are sufficiently smooth. If we work with the ideal space $V_\LOD$ instead of $V_{\LOD}^{\ell,h}$ in \eqref{cG-LOD-scheme}, then the expected convergence of the cG-LOD-method \eqref{cG-LOD-scheme} is for sufficiently smooth solution $u$
\begin{align} \label{convLOD}
	\max_{n = 0,\dots,N} \| u(t_n) - u_\LOD(t_n) \|_{H^1(\D)} \ls \tau^{2q} + H^3.
\end{align}
We note that in the general setting with the fine-scale discretized and localized space $V_{\LOD}^{\ell,h}$, the additional errors are negligible in numerical experiments as long as we choose $h$ sufficiently small and $\ell$ sufficiently large. The optimal convergence rates in \eqref{convLOD} have been verified in numerical experiments (see \cite{DoedingHenningWarnegard}) but we emphasize that a rigorous convergence analysis of the method with sharp regularity assumptions on the solution $u$ is still left for future work. Furthermore, the convergence has not yet been tested for problems with low regularity, i.e., when $V \in L^\infty$ is rough (discontinuous) in the time-dependent GPE. \\
For our numerical comparisons in Section \ref{sec:comp} we define analogously the fully-discrete schemes based on the standard Lagrange FE discretization (cf. \cite{Makridakis,DoedingHenning}): We denote by $\P^k(\T_H)$ the space of standard $H^1$-conforming Lagrange FE of order $k$ associated with the triangulation $\T_H$ and by $P_{H,k}: H^1_0(\D) \rightarrow \P^k(\T_H)$ we denote the corresponding $L^2$-projection. With the spaces
\begin{align*}
	& \tilde{W}^{\ell,h}_{q,k} := \left\{ v : [0,T] \rightarrow \P^k(\T_H) : v_{| I_n} \in \P^q(I_n), \, n = 0,\dots,N-1 \right\}, \\
	& W^{\ell,h}_{q,k} := \tilde{W}^{\ell,h}_{q,k} \cap C([0,T],\P^k(\T_H)),
\end{align*}
the cG-$\P^k$-method seeks $u_{\P^k} \in W^{\ell,h}_{q,k}$ such that $u_{\P^k}(\cdot,0) = P_{H,k}(u_0)$ and
\begin{align} \label{cG-Pk-scheme}
	\int_{0}^T (\i \partial_t u_{\P^k}, v) - (u_{\P^k},v)_\H - \beta \big( |u_{\P^k} |^2 u_{\P^k}, v \big) \, dt = 0 \quad \forall v \in \tilde{W}^{\ell,h}_{q-1,k}
\end{align}
which satisfies for sufficiently smooth solutions $u$ the a-priori estimate (cf. \cite{DoedingHenning})
\begin{align} \label{convPk}
	\max_{n = 0,\dots,N} \| u(t_n) - u_{\P^k}(t_n) \|_{H^1(\D)} \ls \tau^{2q} + H^{k}.
\end{align}
In fact, the solution $u$ needs to be $H^{k+1}$-regular for \eqref{convPk} to hold, see \cite{DoedingHenning} for details.

\sect{Numerical Results}
\label{sec:comp}

\begin{figure}
\centering
\begin{minipage}{\textwidth}
\centering
\begin{minipage}{0.45\textwidth}
\centering
\includegraphics[scale=0.35]{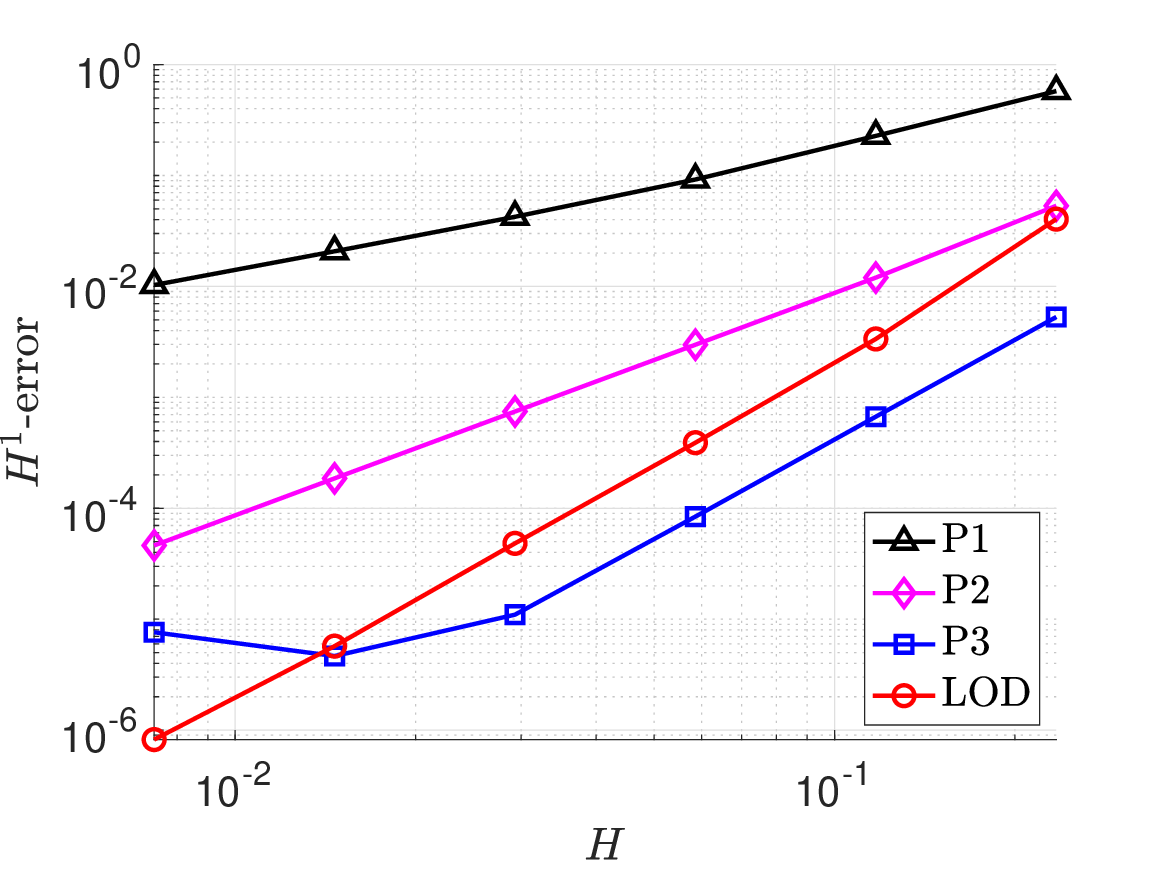}
\caption{Mesh size $H$ vs. $H^1$-error for the continuous potential $V_1$ from \eqref{potentials}.} \label{fig1}
\end{minipage}
\hspace{0.5cm}
\begin{minipage}{0.45\textwidth}
\centering
\includegraphics[scale=0.35]{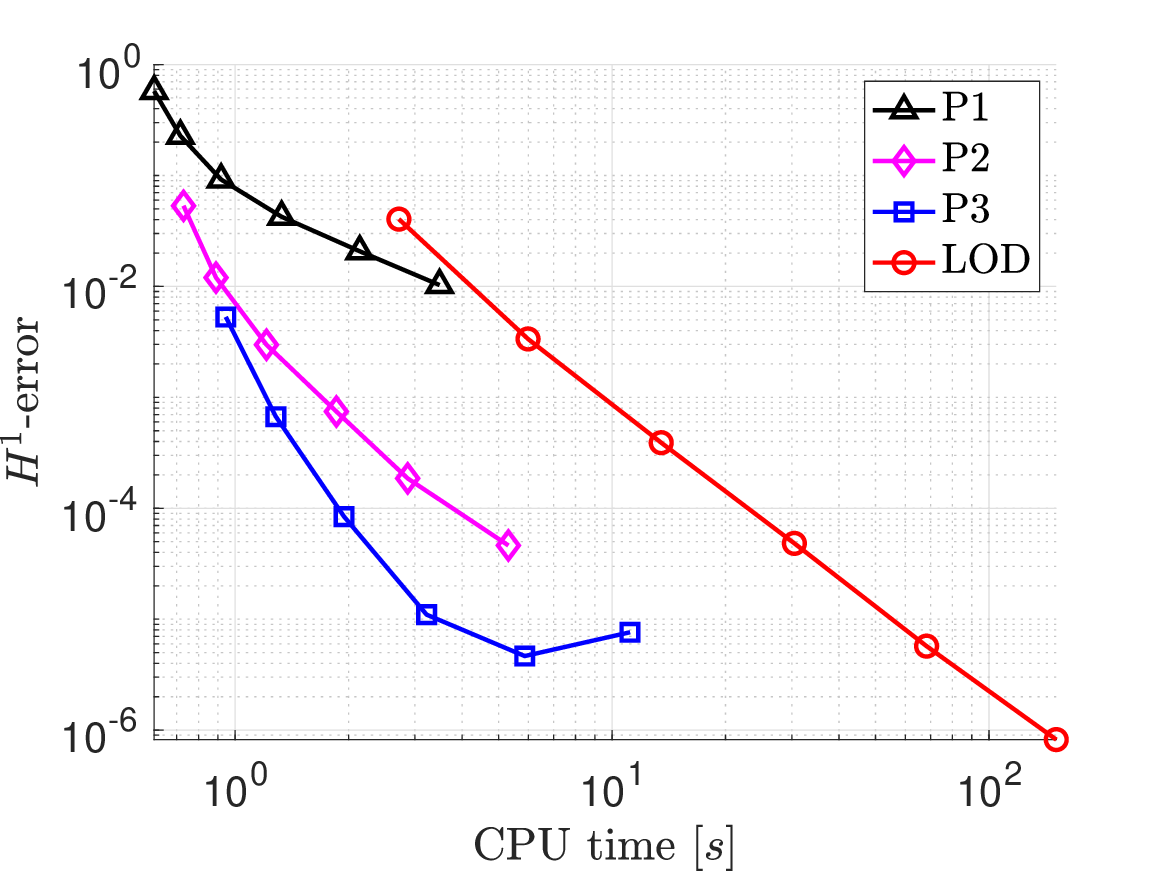}
\caption{CPU time (online) vs. $H^1$-error for the continuous potential $V_1$ from \eqref{potentials}.} \label{fig2}
\end{minipage}
\end{minipage}
\end{figure}

\begin{figure}
\centering
\begin{minipage}{\textwidth}
\centering
\begin{minipage}{0.45\textwidth}
\centering
\includegraphics[scale=0.35]{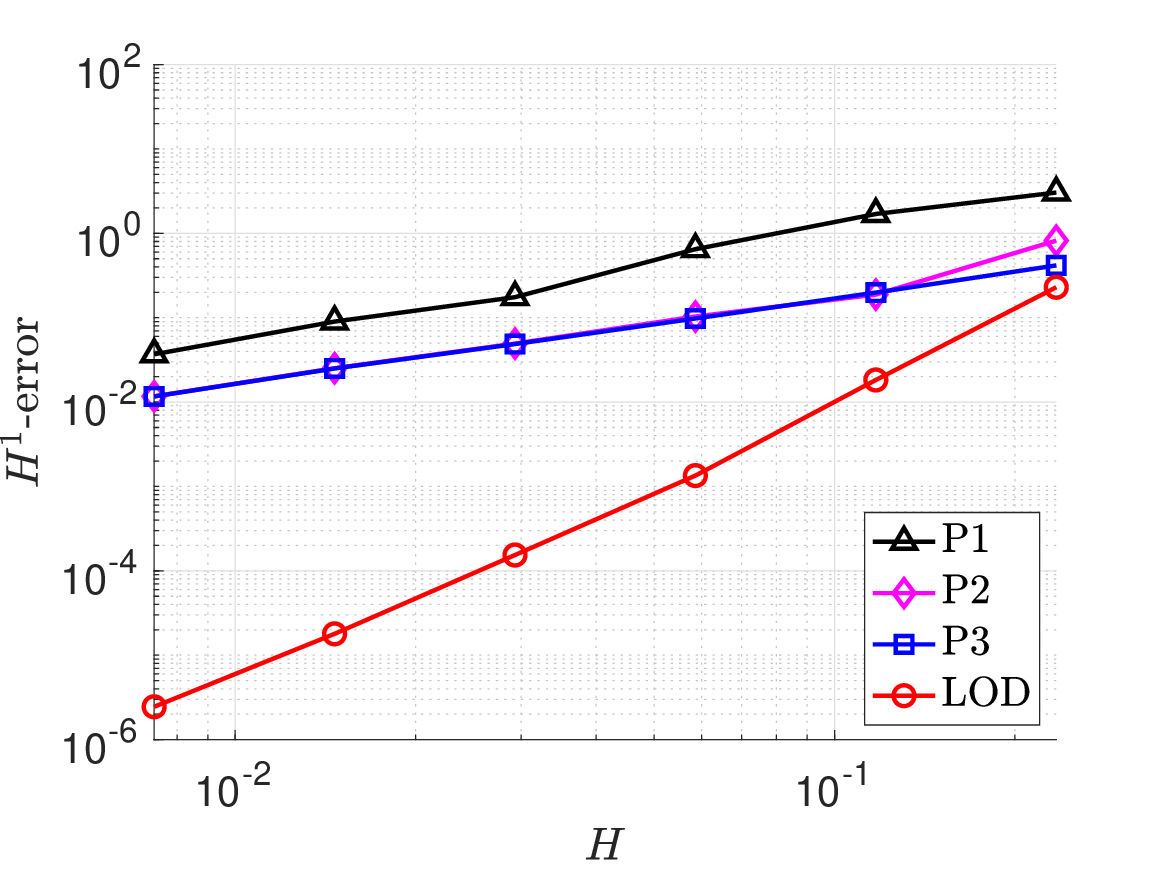}
\caption{Mesh size $H$ vs. $H^1$-error for the discontinuous potential $V_2$ from \eqref{potentials}.} \label{fig3}
\end{minipage}
\hspace{0.5cm}
\begin{minipage}{0.45\textwidth}
\centering
\includegraphics[scale=0.35]{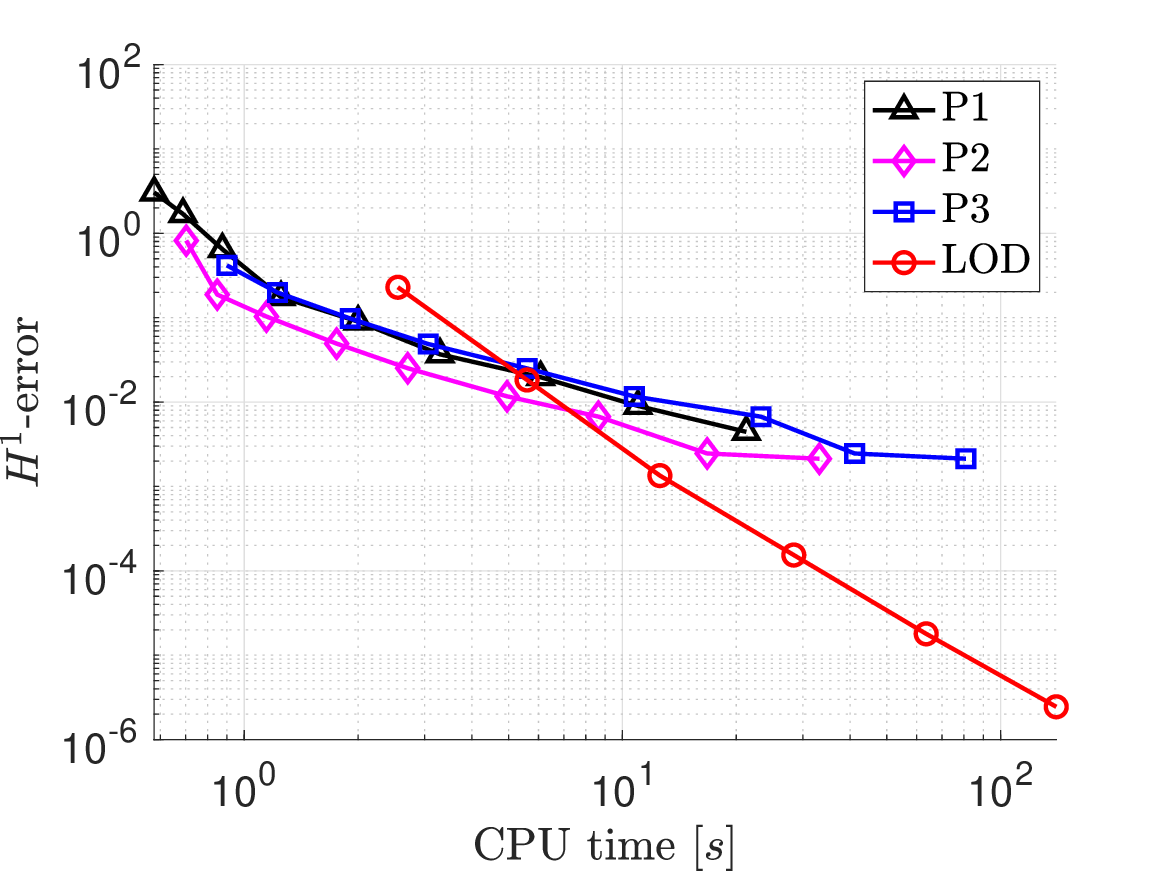}
\caption{CPU time (online) vs. $H^1$-error for the discontinuous potential $V_2$ from \eqref{potentials}.} \label{fig4}
\end{minipage}
\end{minipage}
\end{figure}

In this section we present our numerical results regarding the cG-LOD-method from \eqref{cG-LOD-scheme}, verify its expected convergence from \eqref{convLOD}, and in particular compare its performance with the classical methods from \eqref{cG-Pk-scheme}. All our experiments are run on an Apple MacBook Pro 13" (2022) with an Apple M2 CPU (with 4 $\times $ 3.49 GHz and 4 $\times $ 2.4 GHz cores) and 16GB of physical memory. The implementation is available as a Julia (v1.9.3) code on \url{https://github.com/cdoeding/LODvsFEM}. \\
In the numerical simulations we considers the nonlinear GPE in one space dimension
\begin{align} \label{cont}
	\i \partial_t u = - \partial_{xx} u + V_j u + \beta |u|^2 u \quad \text{in } (-15,15) \times [0,T]
\end{align}
with $\beta = 100$, $T = 0.4$, and two different potentials $V_j$, $j = 1,2$ given by
\begin{align} \label{potentials}
	V_1(x) = 10x^2 \quad \text{and} \quad V_2(x) = (10x^2) \chi_{\{ x \le 0\} }(x) + 100 \chi_{\{ x \ge 5\}}(x) 
\end{align}
where $\chi_{\S}$ denotes the characteristic function on the set $\S \subset \D$. For the initial value we choose the ground state w.r.t. the potential $V_{\gs}(x) = x^2$, i.e., the minimizer of the GPE energy functional $E_{\gs}(v) = \tfrac{1}{2} \int_{\Omega} |\nabla u|^2 + V_{\gs} |u|^2 + \tfrac{\beta}{2} |u|^2 dx$. \\
We solve \eqref{cont} with $V_1$ and $V_2$ from \eqref{potentials} with the cG-LOD-method from \eqref{cG-LOD-scheme} and the cG-$\P^k$-methods from \eqref{cG-Pk-scheme} for $k = 1,2,3$ on uniform meshes. The (coarse) mesh sizes are set to $H = 30 \cdot 2^{-i}$, $i = 7,\dots,12$ and for the cG-LOD-method the fine mesh size is set to $h = 30 \cdot 2^{-16}$ and the oversampling parameter is set to $\ell = i + 5$ for $H = 30 \cdot 2^{-i}$. For all methods we use the continuous Galerkin time integrator with $q = 2$ and step size $\tau = 2 \cdot 10^{-3}$ so that we have a fourth order in time method and the error caused by the time discretization is negligible. According to \cite{DoedingHenning} the nonlinear equations are solved in each time step using a fixed-point iteration with a tolerance of $10^{-10}$ w.r.t. the $L^2$-norm. The initial value $u_0$, given by a minimizer of $E_{\gs}$, is computed by a Sobolev gradient descent method (cf. \cite{BaoCai,HenningPeterseim20}) using a $\P^1$-Lagrange FE discretization on the fine mesh with $h = 30 \cdot 2^{-16}$. \\
For simplicity, we calculate the $H^1$-error via the quantity $\| u_{\mathrm{ref}}(T) - u_{\LOD}(T) \|_{H^1}$ for the cG-LOD-method where $u_{\mathrm{ref}}$ is a reference solution obtained by the cG-$\P^1$-method on the fine mesh with mesh size $h = 30 \cdot 2^{-16}$. The error is calculated accordingly for the cG-$\P^k$-methods with the same reference solution. The CPU times shown in Fig. \ref{fig2} and Fig. \ref{fig4} are in any case measured for the time stepping (online), excluding pre-computations (offline) such as the calculation of the LOD basis functions. There are two reasons for this: First, the calculations of the basis functions are completely independent and can be calculated in parallel. This allows the computational costs to be kept small on modern computer architectures. The second reason is that, especially on long time scales, the costs for the time stepping clearly dominate and the costs for the pre-computations are therefore negligible.

Fig. \ref{fig1} and Fig. \ref{fig2} show the numerical results with the smooth potential $V_1$ from \eqref{potentials}. As expected, we see the convergence behavior of \eqref{convLOD} and \eqref{convPk} as the solution admits $H^4$-regularity for $V \in H^2$. We note that our choice of the oversampling parameter $\ell$ and the fine scale $h$ is sufficient since the error of the cG-LOD-method is not polluted by the localization or fine-scale errors. Only the error of the cG-$\P^3$-method starts to fall off which is due to the choice of our reference solution. Fig. \ref{fig2} visualizes the measured CPU times for the time stepping. We observe that the cG-$\P^3$-method and even the cG-$\P^2$-method are faster than the cG-LOD-method although for fixed $H$, i.e., on a fixed mesh, the dimension of the LOD space is smaller than for the classical methods. This advantage does not pay off because the assembly of the nonlinear terms in the LOD space is too expensive and cancels out this advantage. We believe that this may change in higher dimensions although we have not investigated this further. \\
Fig. \ref{fig3} and Fig. \ref{fig4} show the numerical results with the discontinuous potential $V_2$ from \eqref{potentials}. In this case the solution only admits $H^2$-regularity due to the rough potential and therefore the classical methods show only first order convergence. In contrast, the cG-LOD-method is not affected by the low regularity of the solution since $\H u$ is still $H^2$-regular. Hence the incorporation of the discontinuous potential in the LOD space via 
\eqref{LODspace} preserves the third order of convergence from \eqref{convLOD}. Not surprisingly, due to the higher order, the cG-LOD-method performs better as the classical methods in this case, as shown in Fig. \ref{fig4}. \\
In conclusion, we have shown in numerical experiments that classical methods may be better suited to solve the GPE with smooth coefficients when they show high order convergence. However, for rough problems with non-smooth potentials, we have verified that the cG-LOD-method retains its high order and classical methods cannot compete in these situations. \\

\textbf{Acknowledgement.} This work was funded by the Deutsche Forschungsgemeinschaft (DFG, German Research Foundation) under Germany's Excellence Strategy – EXC-2047/1 – 390685813.

\end{document}